\pdfoutput=1

\documentclass[12pt]{amsart}
\newtheorem{theorem}{Theorem}[section]

\theoremstyle{definition}   
\newtheorem{definition}{Definition}
\newtheorem{example}[theorem]{Example}

\theoremstyle{remark}
\newtheorem{remark}[theorem]{Remark}

\numberwithin{equation}{section}
\usepackage{enumerate}
\usepackage{mathtools}
\usepackage{tikz-cd}
\usepackage{upgreek}
\usetikzlibrary{bending}
\usepackage{graphicx,adjustbox}
\usepackage{hyperref}
\usepackage{amsmath,amssymb}
\usepackage[normalem]{ulem}
\hypersetup{colorlinks,linkcolor={blue},citecolor={red},urlcolor={red}}
\usepackage[a4paper,
            bindingoffset=0.2in,
            left=1in,
            right=1in,
            top=1in,
            bottom=1in,
            footskip=.25in]{geometry}

\usepackage[
backend=biber,
bibstyle=numeric,
sorting=nyt,doi=false,url=false,
]{biblatex}
\DeclareNameAlias{author}{last-first}
\DeclareFieldFormat[article, inbook]{title}{#1} 
\renewbibmacro{in:}{}
\DeclareFieldFormat{journaltitle}{#1\isdot}

\usepackage{amscd}
\usepackage{graphicx}
\usepackage[all]{xy}
\title[$G$-equivariant structures of twisted Higgs bundles]{$G$-equivariant structures of twisted Higgs bundles}
\author{KUNTAL BANERJEE}
\address{Department of Mathematics\\ Swami Vivekananda University Barrackpore India}
\email{kuntalbanerjee53@gmail.com, kuntalb@svu.ac.in}

\begin{document}
\begin{abstract}
In this paper we characterize $G$-equivariant structures of twisted Higgs bundles on $\mathbb{CP}^1$ where $G\subset SL(2, \mathbb{C})$ is a closed reductive subgroup. Also for a reductive group $H\subset PGL(2, \mathbb{C})$ we characterize $H$-equivariant structures of co-Higgs bundles particularly for the case where $H$ does not admit a splitting for the quotient map $SL(2, \mathbb{C})\to PGL(2, \mathbb{C})$.
\end{abstract}
\maketitle 
\tableofcontents

\let\thefootnote\relax\footnotetext{2020\textit{ Mathematics Subject Classification.} 14H60, 14F06,  14L30.}\let\thefootnote\relax\footnotetext{}\let\thefootnote\relax\footnotetext{\textit{Keywords and phrases.} Complex projective line, Harder-Narasimhan filtration, Haar measure}

\section{Introduction}

\noindent Let $G$ be a reductive closed subgroup of $SL(2, \mathbb{C})$. In this paper we define a $G$-invariant structure of $\mathcal{O}_{\mathbb{CP}^1}(d)$-twisted Higgs bundles on $\mathbb{CP}^1$ extending similar definitions from \cite{BKP}. Since $G$ is a reductive group $G$ admits a maximal compact subgroup $K_G$ with a Haar measure $\mu_{K_G}$. This Haar measure defines a $G$-equivariant holomorphic splitting for every $\mathcal{O}_{\mathbb{CP}^1}(d)$-twisted Higgs bundle $\mathcal{E}$ even if $\mathcal{E}$ does not have a $G$-equivariant structure. Now consider a $G$-equivariant twisted Higgs bundle $\mathcal{E}$. Then $\mathcal{E}$ is $G$-equivariantly holomorphically isomorphic to a direct sum $\oplus_{i=1}^{k}\oplus_{j=1}^{r_i} \mathcal{O}(a_i)\otimes_\mathbb{C} M_i, [\Phi_{ij}]$ where $\Phi_{ij}\in \left(H^0(\mathbb{P}^1, \mathcal{O}(a_j - a_i + d))\otimes Hom(M_i, M_j)\right)^G$ where $M_i$ is a finite dimensional complex $G$-module $M_i$ (Theorem \ref{maintheorem2}). Similar characterization results are obtained for $H$-equivariant co-Higgs bundles where $H\subset PGL(2, \mathbb{C})$ for two cases based on splitting of $H$ for the quotient map $Aut(\mathcal{O}_{\mathbb{CP}^1}(1)) = SL(2, \mathbb{C})\to PGL(2, \mathbb{C})$ (Theorem \ref{maintheorem3}, Theorem \ref{maintheorem4}).\\

\noindent A co-Higgs bundle $\mathcal{E} = (E, \Phi)$ on $\mathbb{CP}^1$ is a holomorphic vector bundle $E\to \mathbb{CP}^1$ and a holomorphic bundle morphism $E\to E\otimes \mathcal{O}_{\mathbb{CP}^1}(2)$. For a detailed analysis for geometry of co-Higgs bundles refer to \cite{Ste1} and \cite{Rayan}. If $H$ admits a splitting then $H$-equivariant co-Higgs bundles can be characterized up to an $H$-equivariant isomorphism similarly as Theorem \ref{maintheorem2}. Now if $H$ does not split in to $SL(2, \mathbb{C})$ the line bundle $\mathcal{O}_{\mathbb{CP}^1}(2)$ admits an $H$-equivariant structure. This structure is essentially used in defining $H$-equivariant structures on co-Higgs bundles $\mathcal{E}$. In this case $-I_2\in \beta^{-1}(H)$ acts strictly via multiplication by $-1$ on $\mathcal{O}_{\mathbb{CP}^1}(1)$ where $\beta: SL(2, \mathbb{C})\to PGL(2, \mathbb{C})$ is the quotient morphism. An $H$ equivariant co-Higgs bundle which is of the form $(\oplus_{i=1}^r \mathcal{O}(a)\otimes_{\mathbb{C}}M, \Phi)$ while $a$ is an even number is given by an $H$-equivariant module $M$ and $\Phi\in (H^0(\mathcal{O}(2))\otimes End(M))^H$. On the other hand if $a$ is odd then $\mathcal{O}(a)$ has a $\beta^{-1}(H)$-module structure with $-I_2$ acting via multiplication of $-1$ and for a finite dimensional complex $\beta^{-1}(H)$-module $M$ with $-I_2$ acting via multiplication of $-1$ there is an $H$-equivariant structure on $\mathcal{O}(a)\otimes M$. Similarly there is an $H$-equivariant structure on $H^0(\mathcal{O}(2))\otimes_{\mathbb{C}} End(M)$ which contains $H$-equivariant co-Higgs fields $\Phi$ on $\mathcal{O}(a)\otimes M$. From these observations theorems \ref{maintheorem3} and \ref{maintheorem4} follow.

\section{Equivariant twisted Higgs bundles}
\noindent Consider the standard $SL(2, \mathbb{C})$-action $\hat{\rho}$ on $\mathbb{CP}^1$ induced by the natural action $\rho: SL(2, \mathbb{C})\times\mathbb{C}^2\to\mathbb{C}^2$. In other words $\hat{\rho}\in Aut(\mathbb{CP}^1)$. Then there exists an isomorphism $\hat{\sigma}:PSL(2, \mathbb{C})\to Aut(\mathbb{CP}^1)$. Now let $G\subset SL(2, \mathbb{C})$ be a closed reductive subgroup.
\begin{definition}\label{equiv}
A $G$-equivariant holomorphic vector bundle is a holomorphic vector bundle $E\to \mathbb{P}^1$ and a $\rho$-compatible action $\phi$ which satisfies the following commutative diagram with respect to the bundle projection of $E$
\begin{equation}
\begin{tikzcd}
G\times E \arrow{r}{\phi} \arrow[swap]{d} & E \arrow{d}{} \\
G\times\mathbb{CP}^1 \arrow{r}{\rho} & \mathbb{CP}^1
\end{tikzcd}
\end{equation}

\noindent and the map $E\to E$ defined by $v\mapsto \phi(g, v)$ is an isomorphism of vector bundles over $\hat{\rho}(g)$, that is, the following diagram is commutative. Now on wards, we will use the symbol $g\cdot e$ alternatively to $\phi(g, e)$.

\begin{equation}
\begin{tikzcd}
E \arrow{r}{\phi(g, \cdot)} \arrow[swap]{d} & E \arrow{d}{} \\
\mathbb{CP}^1 \arrow{r}{\hat{\rho}(g)} & \mathbb{CP}^1
\end{tikzcd}
\end{equation}
\end{definition}
\noindent The action of $SL(2, \mathbb{C})$ on $\mathbb{CP}^1$ induces an action on the hyperplane line bundle $\mathcal{O}_{\mathbb{CP}^1}(1)$ in the following way. Note that every element $g\in SL(2, \mathbb{C})$ acts on $\mathbb{CP}^1$ in terms of M\"obius transformations. If $[z_0, z_1]\in \mathbb{CP}^1$ and $g = \begin{bmatrix}
    a & b\\
    c & d
\end{bmatrix}\in SL(2, \mathbb{C})$ then $$g\cdot [z_0, z_1] = [az_0 + bz_1, cz_0 + dz_1].$$ Then on the tautological line bundle $\mathcal{O}_{\mathbb{CP}^1}(-1)$ there is a group action which on each fiber reads as $g\cdot ([z_0, z_1], v) = (g\cdot [z_0, z_1], g(v))$. This action induces the dual action on $\mathcal{O}_{\mathbb{CP}^1}(1)$ by $g\cdot ([z_0, z_1], f) = (g\cdot [z_0, z_1], f\circ g^{-1})$. These group actions define $G$-equivariant structures on every line bundle $\mathcal{O}_{\mathbb{CP}^1}(d)$ by tensoring the action $d$ times where $d \in\mathbb{Z}$.

\begin{definition}\label{equiv1}
Let $\mathcal{E} = (E, \Phi)$ be any $\mathcal{O}_{\mathbb{CP}^1}(d)$- twisted Higgs bundle where $d\geq 2$ where $E$ is a $G$-equivariant vector bundle. Note that $G$ acts on $E(d) := E\otimes\mathcal{O}_{\mathbb{CP}^1}(d)$ by $(\phi\otimes\rho)(e\otimes l) = \phi(e)\otimes\rho(l)$ where $e\otimes l\in E(d)$. Then $\mathcal{E}$ is said to be $G$-equivariant with respect to the natural $G$-equivariant structure of $\mathcal{O}(d)$ if $\Phi$ is $\rho$-compatible, that is, $\Phi$ satisfies the following commutative diagrams.

\begin{equation}
\begin{tikzcd}[row sep=2em, column sep=2em]
  & G\times E \arrow[rr, "\phi"] \arrow[dd, "1_G\times \text{proj}" near start] \arrow[dl, "1_G\times\Phi"'] & & E \arrow[dd, "\text{proj}"] \arrow[dl, "\Phi"] \\
  G\times E(d) \arrow[rr, crossing over, "(\phi\otimes\rho)" near start] \arrow[dd, "1_G\times\text{proj}"'] & & E(d) \\
  & G\times\mathbb{CP}^1 \arrow[rr, "\rho" near start] \arrow[dl, "1_G\times 1_{\mathbb{CP}^1}"'] & & \mathbb{CP}^1 \arrow[dl, "1_{\mathbb{CP}^1}"] \\
  G\times\mathbb{CP}^1 \arrow[rr, "\rho"'] & & \mathbb{CP}^1 \arrow[from=uu, crossing over, "\text{proj}" near start]
\end{tikzcd}
\end{equation}
\end{definition}

\begin{definition}\label{equiv2}
If $\mathcal{E} = (E,\phi_E, \Phi)$ and $\mathcal{F} = (F, \phi_F, \Psi)$ are two $G$-equivariant $\mathcal{O}_{\mathbb{CP}^1}(d)$-twisted Higgs bundles then $\mathcal{E}$ and $\mathcal{F}$ are isomorphic as $G$-equivariant twisted Higgs bundles, also written as, $\mathcal{E}\cong_G\mathcal{F}$ if the following diagrams commute for some isomorphism $f:E\to F$ given as follows.

\begin{equation}
\begin{tikzcd}
G\times E \arrow{r}{\phi_E} \arrow[swap]{d}{1_G\times f} & E  \arrow{r}{\Phi} \arrow[swap]{d}{f} & E(d) \arrow[swap]{d}{f\otimes 1_{\mathcal{O}_{\mathbb{CP}^1}}(d)}\\
G\times F \arrow{r}{\phi_F} & F \arrow{r}{\Psi} & F(d)
\end{tikzcd}
\end{equation}
\end{definition}

\begin{example}\label{equivariant}
Let $a\in \mathbb{Z}$ and $\mathcal{E} = \bigoplus\limits_{i=1}^n (\mathcal{O}(a), \Phi_i)$ where $\Phi_i\in H^0(\mathbb{CP}^1, \mathcal{O}(d))$. Consider the $G$-equivariant structure $\lambda$ on $\mathcal{O}_{\mathbb{CP}^1}(d)$, which is a family of isomorphisms $\lambda_{g, x}: \mathcal{O}(d)_x\to \mathcal{O}(d)_{g\cdot x}$ for $x\in\mathbb{P}^1$ and $g\in G$ that follows $\lambda_{gh, x} = \lambda_{g, h\cdot x}\circ \lambda_{h, x}$. Then $s\in H^0(\mathbb{CP}^1, \mathcal{O}(d))$ is acted by $\lambda$ where the action is given by $$(g\cdot s)(x) = \lambda_{g, g^{-1}x}(s(g^{-1}x)).$$ Thus 
$\mathcal{E}$ is $G$-equivariant if and only if $\Phi_i$ is invariant under the action of $G$.
\end{example}

\begin{theorem}
Let $\mathcal{E}$ be a $G$-equivariant twisted Higgs bundle. Let $0 = \mathcal{E}_0\subsetneq \mathcal{E}_1\subsetneq\dots\subsetneq \mathcal{E}_{l-1}\subsetneq \mathcal{E}_l = \mathcal{E}$. Then the action of $G$ preserves the Harder-Narasimhan filtration.   
\end{theorem}
\begin{proof}
Note that $0 = \phi(g, \mathcal{E}_0)\subset\dots \phi(g, \mathcal{E}_{l-1})\subset \phi(g, \mathcal{E}_l) \cong \mathcal{E}$ is a Harder-Narasimhan filtration of $\mathcal{E}$ for any Harder-Narasimhan filtration $0 = \mathcal{E}_0\subset\dots \mathcal{E}_{l-1}\subset \mathcal{E}_l \cong \mathcal{E}$ the following filtration satisfies all the requirements of a Harder-Narasimhan filtration
\begin{equation}
0 = \phi(g, E_0)\subset\dots \phi(g, E_{l-1})\subset \phi(g, E_l) \cong \mathcal{E}
\end{equation}
since $\Phi(\phi(g, E_i))\subset(\phi\otimes\rho)(g, E_i\otimes\mathcal{O}(d))$. By the uniqueness of the Harder-Narasimhan filtration the statement follows.
\end{proof}

\begin{theorem}
Let $0\to \mathcal{E}'\to\mathcal{E}\to \mathcal{E}''\to 0$ be a short exact sequence of $G$-equivariant twisted Higgs bundles on $\mathbb{CP}^1$. If this exact sequence is a holomorphic splitting of twisted Higgs bundles then it induces a $G$-equivariant splitting, that is, $\mathcal{E}\cong \mathcal{F}\oplus\mathcal{G}$ for some $G$-equivariant twisted Higgs bundles $\mathcal{F}$ and $\mathcal{G}$. 
\end{theorem}
\begin{proof}
The proof is similar to the proof of Proposition 2.3 of \cite{BKP}. Let $\psi:\mathcal{E}''\to \mathcal{E}$ be a splitting, that is, section of the projection $q:\mathcal{E}\to \mathcal{E}''$. Consider the homomorphism $\psi_g: \mathcal{E}''\to \mathcal{E}$ by $v\mapsto g^{-1}\psi(g\cdot v)$. Thus at the level of bundles $q:E\to E''$ satisfies the commutativity relation $(q\otimes 1_{\mathcal{O}_{\mathbb{CP}^1}(d)})\circ \Phi = \Phi_{E''}\circ q$. Also the homomorphism $\psi$ satisfies $\Phi\circ \psi = (\psi\otimes 1_{\mathcal{O}_{\mathbb{CP}^1}(d)})$. Now consider a compact subgroup $K_G$ contained in $G$ and a Haar measure $\mu$ on $K_G$. Define the bundle homomorphism 

$$\tilde{\psi}(v) = \int_{K_G} \psi_g(v)d\mu_{K_G}(g) = \int_{K_G} (g^{-1} \cdot \psi(g\cdot v))d\mu_{K_G}(g).$$ 

\noindent Fiberwise on $E''$ we must show that $\Phi(\tilde{\psi}(v)) = (\tilde{\psi}\otimes 1_{\mathcal{O}_{\mathbb{CP}^1}(d)})(\Phi_{E''}(v))$. At first note that $$\Phi(g^{-1}\cdot\psi(g\cdot v)) = (g^{-1}\otimes \rho(g^{-1}))\cdot \Phi(\psi(g\cdot v)).$$ Since $\Phi\circ \psi = (\psi\otimes 1_{\mathcal{O}_{\mathbb{CP}^1}(d)})\circ \Phi_{E''}$ we have $$(g^{-1}\otimes \rho(g^{-1}))\cdot \Phi(\psi(g\cdot v)) = (g^{-1}\otimes \rho(g^{-1}))\cdot (\psi\otimes 1_{\mathcal{O}_{\mathbb{CP}^1}(d)})( \Phi_{E''}(g\cdot v)).$$ 

\noindent Now $\Phi_{E''}(g\cdot v) = (g\otimes\rho(g))(\Phi_{E''}(v))$ and therefore $$\Phi(\psi_g(v)) = (g^{-1}\otimes \rho(g^{-1}))\cdot (\psi\times 1_{\mathcal{O}_{\mathbb{CP}^1}(d)}) (g\otimes \rho_g)(\Phi_{E''}(v)).$$ On the $\mathcal{O}_{\mathbb{CP}^1}(d)$ twist $\rho_{g^{-1}}\circ 1_{\mathcal{O}_{\mathbb{CP}^1}(d)}\circ \rho_g = 1_{\mathcal{O}_{\mathbb{CP}^1}(d)}$ and $g^{-1}\circ \psi\circ g = \psi_g$. Now 
\begin{align*}
&\Phi(\tilde{\psi}(v)) = \Phi(\int_{K_g} \psi_g(v) d\mu_{K_G}(g)) \\
&= \int_{K_G} \Phi(\psi_g(v))d\mu_{K_G}(g) \\
&= \int_{K_G}(\psi_g\otimes 1_{\mathcal{O}_{\mathbb{CP}^1}(d)})\Phi_{E''}(v) d\mu_{K_G}(g) 
\end{align*}

\noindent and
\begin{align*}
&\int_{K_G}(\psi_g\otimes 1_{\mathcal{O}_{\mathbb{CP}^1}(d)})\Phi_{E''}(v) d\mu_{K_G}(g) \\
&= (\int_{K_G}\psi_g d\mu_{K_G}(g)\otimes 1_{\mathcal{O}_{\mathbb{CP}^1}(d)})(\Phi_{E''}(v) \\
&= (\tilde{\psi}\otimes 1_{\mathcal{O}_{\mathbb{CP}^1}(d)})(\Phi_{E''}(v)).
\end{align*}
\noindent It requires to verify that $\tilde{\psi}$ is a splitting. Fiberwise at each $x\in \mathbb{P}^1$ we have $$q_x(\psi_g(v)) = q_x(g^{-1}\cdot\psi_{g\cdot x}(g\cdot v)) = g^{-1} q_{g\cdot x}(\psi_{g\cdot x}(g\cdot v)) = g^{-1}\cdot (g\cdot v) = v.$$ Integrating over $K_G$ we have $$(q\circ \tilde{\psi})(v) = q_x(\int_{K_G} \psi_g(v)d\mu_{K_G}(g)) = \int_{K_G} q_x(\psi_g(v))d\mu_{K_G}(g) = \int_{K_G} v d\mu_{K_G} = v.$$ Hence the proof of the theorem.
\end{proof}

\begin{remark}
Let $\mathcal{E} = \bigoplus\limits_{i=1}^k \mathcal{E}_i$ be a $G$-equivariant twisted Higgs bundle where $\mathcal{E}_i = \left(\bigoplus\limits_{j=1}^{r_i} \mathcal{O}(a_i), [\Phi_{ij}]\right)$ and $\Phi_{ij}$ is a global matrix representation of $\mathcal{O}_{\mathbb{CP}^1}(d)$-twisted Higgs field on $r_i$-fold summands $$\mathcal{O}(a_i)\to \mathcal{O}(a_i)\otimes\mathcal{O}_{\mathbb{CP}^1}(d).$$

\noindent Moreover assume that $a_1 >\dots > a_k$. Then $\mathcal{E}$ admits a $G$-equivariant Harder-Narasimhan filtration of twisted Higgs bundles which is given as $$0\subset \mathcal{E}_k\subset \mathcal{E}_k\oplus\mathcal{E}_{k-1}\subset\dots\subset \mathcal{E}_k\oplus\dots\oplus\mathcal{E}_1 = \mathcal{E}.$$ Then $\mathcal{E}$ admits a $G$-equivariant splitting $\mathcal{E} \cong \bigoplus_{i=1}^k \mathcal{E}_i'$. In particular, $\mathcal{E}_i\cong_G \mathcal{E}_i'$ if and only if $\mathcal{E}_i$ is a $G$-equivariant $\mathcal{O}_{\mathbb{CP}^1}(d)$-twisted Higgs bundle in Example \ref{equivariant}. 
\end{remark}

\begin{theorem}\label{maintheorem1}
Let $\mathcal{E} = (\oplus_{i=1}^r \mathcal{O}(a), \Phi)$ be an $\mathcal{O}_{\mathbb{CP}^1}(d)$-twisted $G$-equivariant Higgs bundle with induced $G$-equivariant structure on $\oplus_{i=1}^r \mathcal{O}(a)$. Then there is a finite dimensional complex $G$-module $M$ such that $\mathcal{E}$ is $G$-equivariantly holomorphically isomorphic to $(\mathcal{O}(a)\otimes_{\mathbb{C}} M, \Phi_0)$ where $$\Phi_0\in \left(H^0(\mathbb{CP}^1, \mathcal{O}_{\mathbb{CP}^1}(d))\otimes_{\mathbb{C}} End(M)\right)^G$$ and $V^G$ is the space of $G$-invariant elements for a complex $G$-module $V$. Moreover, $M$ is unique up to a $G$-equivariant isomorphism of complex vector spaces.\\

\noindent Let $\mathcal{E} = \oplus_{i=1}^n \mathcal{E}_i$ where each $\mathcal{E}_i$ is a $G$-equivariant $\mathcal{O}_{\mathbb{CP}^1}(d)$-twisted Higgs bundle of the type mentioned in the above paragraph and $a_i\neq a_j$ for all $i\neq j$ where $\mathcal{O}(a_i)$ is the underlying line bundle of $\mathcal{E}_i$. Then $a_i$'s are unique integers and $M_i$'s are unique complex $G$-modules.
\end{theorem}

\begin{proof}
Let $W = \oplus_{i=1}^r \mathcal{O}(a)$ and $M = H^0(\mathbb{CP}^1, Hom(\mathcal{O}(a), W))$. Here $M$ is also the complex vector space of global bundle homomorphisms $\mathcal{O}(a)\to W$. Here $M$ is a $G$-module with an induced action from respective $G$-actions on $\mathcal{O}(a)$ and $W$. The $G$-equivariant homomorphism $f: \mathcal{O}(a)\otimes M\to W$ is fiberwise defined by the evaluation homomorphism $w\otimes m\to m(w)\in W$. Also consider an element $\Phi_0\in \left(H^0(\mathbb{CP}^1, \mathcal{O}_{\mathbb{CP}^1}(d))\otimes_{\mathbb{C}} End(M)\right)^G$ which is of the form 
\begin{equation}\label{equivariant Higgs fields}
\Phi_0 = \sum_{k=1}^{r^2} \Phi_0^{(k)}\otimes \psi_0^{(k)}. 
\end{equation}
Then for each value of $k$ the element $\Phi_0^{(k)}\otimes \psi_0^{(k)}$ defines a $G$-equivariant $\mathcal{O}(d)$-twisted Higgs field $\mathcal{O}(a)\otimes M\to \mathcal{O}(a+d)\otimes M$, so defines $\Phi_0$, say, $f_{\Phi_0}$. The evaluation homomorphism on $f_{\Phi_0}$ then defines an $\mathcal{O}(d)$-twisted $G$ equivariant Higgs bundle $$\Phi: W\to W\otimes\mathcal{O}(d).$$ 

\noindent Let $\mathcal{O}(a)\otimes M\cong_G \mathcal{O}(a)\otimes M'$ for complex $G$ modules $M$ and $M'$. If the above isomorphism is restricted to point wise fibers of $\mathcal{O}(a)$ there is a $G$-equivariant isomorphism of complex vector spaces $M$ and $M'$. Uniqueness of the integers $a_i$'s and $M_i$'s follow from the respective Harder-Narasimhan filtration of $\mathcal{E}$.
\end{proof}
\begin{remark}\label{directsum}
For arbitrary integers $a \neq b$ consider two vector bundles $E=\oplus_{i=1}^m \mathcal{O}(a)$ and $F=\oplus_{j=1}^n \mathcal{O}(b)$. Then on $E\oplus F$ the $G$-equivariant structure is defined by the respective structures of $\mathcal{O}(a)\otimes_{\mathbb{C}} M$ and $\mathcal{O}(b)\otimes_{\mathbb{C}} N$ where $M$ and $N$ are finite dimensional complex $G$-modules. Then the cross diagonal components of $G$-equivariant $\mathcal{O}_{\mathbb{CP}^1}$-twisted Higgs bundles are given by elements in $$\left(H^0(\mathcal{O}(a-b+d))\otimes_\mathbb{C} Hom_\mathbb{C}(N, M)\right)^G~\text{and}~ \left(H^0(\mathcal{O}(b-a+d))\otimes_\mathbb{C} Hom_\mathbb{C}(M, N)\right)^G~$$ where $Hom_\mathbb{C}(N, M)$ and $Hom_\mathbb{C}(M, N)$ have the canonical complex $G$-module structures coming from the respective $G$-module structures of $M$ and $N$. 
\end{remark}
\begin{theorem}\label{maintheorem2}
Every $G$-equivariant $\mathcal{O}(d)$-twisted Higgs bundle $\mathcal{E}$ is $G$-equivariantly and holomorphically isomorphic to a direct sum $\oplus_{i=1}^{k}\oplus_{j=1}^{r_i} \mathcal{O}(a_i)\otimes_\mathbb{C} M_i, [\Phi_{ij}]$ where $\Phi_{ij}\in \left(H^0(\mathbb{P}^1, \mathcal{O}(a_j - a_i + d))\otimes Hom(M_i, M_j)\right)^G$ where $M_i$ is a finite dimensional complex $G$-module $M_i$.
\end{theorem}
\begin{proof}
The proof follows by combining the Birkhoff- Grothendieck decomposition of holomorphic vector bundles on $\mathbb{CP}^1$ with Theorem \ref{maintheorem1} and Remark \ref{directsum}.
\end{proof}
\section{Equivariant co-Higgs bundles under action of $H\subset PGL(2, \mathbb{C})$}

\noindent Let $H\subset PGL(2, \mathbb{C})$ be a reductive subgroup. An $H$ equivariant co-Higgs bundle $\mathcal{E}$ is defined by combining definitions \ref{equiv}, \ref{equiv1} and \ref{equiv2} after replacing the group $G$ is by a reductive subgroup $H\subset PGL(2, \mathbb{C})$. The standard action $SL(2, \mathbb{C})\times\mathbb{C}^2\to\mathbb{C}^2$ produces an action $PGL(2, \mathbb{C})\times \mathcal{O}_{\mathbb{CP}^1}(2)\to \mathcal{O}_{\mathbb{CP}^1}(2)$ since the center $\{\pm I_2\}\subset SL(2, \mathbb{C})$ acts trivially on $\mathcal{O}_{\mathbb{CP}^1}(2)$. In fact the fiberwise action is given by $$g\cdot([z_0, z_1], f\otimes f') = (g\cdot [z_0, z_1], (f\circ g^{-1})\otimes (f'\circ g^{-1})).$$ 

\noindent Let $0\to \{\pm I_2\}\to SL(2,\mathbb{C})\xrightarrow[]{\beta} SL(2, \mathbb{C})/\{\pm I_2\} = PGL(2, \mathbb{C})$ be the exact sequence of the determinant morphism. Then let us assume at first that the exact sequence $0\to \{\pm I_2\}\to \beta^{-1}(H)\to H\to 0$ splits, that is, there exists a homomorphism $\gamma:H\to\beta^{-1}(H)$ such that the composition map $H\xrightarrow[]{\gamma} \beta^{-1}H\xrightarrow[]{\beta} H$ is the identity map. It is proved that $\mathcal{O}_{\mathbb{CP}^1}(1)$ has an $H$-equivariant structure precisely if $H$ has a splitting. If $H$ has a splitting then $\mathcal{O}_{\mathbb{CP}^1}(1)$ has a $\gamma(H)$-equivariant action thus an $H$-equivariant action. On the other hand let $\mathcal{O}_{\mathbb{CP}^1}(1)$ be with an $H$-equivariant structure thus there is a homomorphism $\alpha:H\to Aut(\mathcal{O}_{\mathbb{CP}^1}(1)) = SL(2, \mathbb{C})$ which defines a splitting of $H$. This equivalent criteria of existence of $G$-equivariant structures on $\mathcal{O}(2a + 1)$ can be proved for any integer $a$. Also, the action of $H$ leaves the corresponding Harder-Narasimhan filtration invariant and every bundle $E$ admits a splitting into $H$-equivariant vector bundles. Similar constructions of $H$-equivariant structures can be produced for $\mathcal{O}_{\mathbb{CP}^1}(d)$-twisted Higgs bundles of the type $(\oplus_{i=1}^r \mathcal{O}(a), \Phi)$
if $H$ has a splitting. If $H$ does not have a splitting then note that $\mathcal{O}(a)$ and $\mathcal{O}(a+2)$ either have $H$-equivariant structure or do not have any $H$-equivariant structure together. So we should be careful in defining $H$-equivariant structures on $\mathcal{E}$ when $a$ is odd. We will use the fact that if $E$ has a $H$ module structure then it has a $\beta^{-1}(H)$ module structure with $-I_2$ acting trivially on $E$ and vice versa.

\begin{theorem}\label{maintheorem3}
Suppose that $H$ does not have a splitting $\gamma$ and $a$ be an odd integer then every $H$-equivariant co-Higgs bundle $\mathcal{E} = (\oplus_{i=1}^r \mathcal{O}(a), \Phi)$ is $H$-equivariantly and holomorphically isomorphic to $(\mathcal{O}(a)\otimes_\mathbb{C} M, \Phi_0)$ where $M$ is a complex $\beta^{-1}(H)$-module unique up to isomorphism and $-I_2\in H$ acts on $M$ by multiplication of $-1$ and $\Phi_0\in (H^0(\mathcal{O}(2))\otimes_{\mathbb{C}} End(M))^H$.
\end{theorem}

\begin{proof}
Let $M = H^0(\mathbb{P}^1, Hom(\mathcal{O}(a), W))$ where $W = \oplus_{i=1}^r\mathcal{O}(a)$. The $H$-equivariant structure of $\mathcal{O}(2)$ is the structure defined by $PGL(2, \mathbb{C})$. The $\beta^{-1}(H)$-equivariant action on $\mathcal{O}(a)$ and $\mathcal{O}(a+2)$ is given by the definition that $-I_2$ acts by action of $-1$. Take the $\beta^{-1}(H)$-module $M = H^0(\mathbb{P}^1, Hom(\mathcal{O}(a), W))$ where $W = \oplus_{i=1}^r \mathcal{O}(a)$ has induced $\beta^{-1}(H)$ action and $-I_2$ acts by multiplication of $-1$ on $M$. Also the group $\beta^{-1}(H)$ acts on $End(M)$ where $-I_2$ acts by multiplication of $1$ on $M$. Thus, under these group actions there is a unique $H$-equivariant structure on $\mathcal{O}(a)\otimes_\mathbb{C} M$. Now the evaluation morphism $\mathcal{O}(a)\otimes_{\mathbb{C}} M\to W$ sets a $\beta^{-1}(H)$-equivariant isomorphism. But $-I_2$ acts via the multiplication by $1$ on $W$ since $W$ has an $H$-equivariant structure. So the evaluation isomorphism is $H$-equivariant. Similarly, there is a unique $H$-equivariant structure on $W\otimes\mathcal{O}(2)$ with the $G$-equivariant evaluation isomorphism $\mathcal{O}(a+2)\otimes_{\mathbb{C}} M\to W(2)$.\\

\noindent Let us now consider the construction of $H$-equivariant structure on a co-Higgs field $W\to W(2)$ defined by a $G$-equivariant co-Higgs field $\mathcal{O}(a)\otimes_{\mathbb{C}} M\to \mathcal{O}(a+2)\otimes_{\mathbb{C}} M$ as in Equation \ref{equivariant Higgs fields}. Then each $G$-equivariant co-Higgs field $\Phi_0 :\mathcal{O}(a)\otimes_{\mathbb{C}} M\to \mathcal{O}(a+2)\otimes_{\mathbb{C}} M$ is an element inside $(H^0(\mathbb{CP}^1, \mathcal{O}_{\mathbb{CP}^1}(2))\otimes_{\mathbb{C}} End(M))^H$. The proof of uniqueness of $M$ follows from the same argument of uniqueness used in Theorem \ref{maintheorem1}.
\end{proof}

\begin{theorem}\label{maintheorem4}
(i) Suppose that $H$ does not have a splitting $\gamma$ and $a$ be an odd integer and $b$ be an even integer. Then a $G$-equivariant bundle homomorphism $\Phi:\mathcal{O}(a)\otimes M\to\mathcal{O}(b+2)\otimes N$ comes from $$\left(H^0(\mathcal{O}(b-a+2))\otimes_\mathbb{C} Hom_\mathbb{C}(M, N)\right)^H$$ where $M$ is acted by $-I_2$ via multiplication of $-1$ and $N$ is acted by the multiplication of $1$. \\

\noindent (ii) Let $\mathcal{E}$ be an $H$-equivariant Higgs bundle. Then there is a set of odd integers and $G$-modules $\{a_i, M_i\}_{i=1}^m$ and there is a set of even integers and $G$-modules $\{b_j, N_j\}_{j=1}^n$ such that $-I_2$ acts by multiplication of $-1$ on each $\mathcal{O}(a_i)$ and on each $M_i$, and $-I_2$ acts by multiplication by $1$ on each $\mathcal{O}(b_j)$ and on each $N_j$. Also the homomorphism blocks of $\Phi$ are given by 

$$\begin{cases}
\Phi_{ii'}\in \left(H^0(\mathcal{O}(a_i-a_{i'}+2))\otimes_\mathbb{C} Hom_\mathbb{C}(N_j, M_i)\right)^H\\
\Phi_{jj'}\in \left(H^0(\mathcal{O}(b_j-b_{j'}+2))\otimes_\mathbb{C} Hom_\mathbb{C}(N_j, M_i)\right)^H\\
\Phi_{ij}\in \left(H^0(\mathcal{O}(a_i-b_j+2))\otimes_\mathbb{C} Hom_\mathbb{C}(N_j, M_i)\right)^H\\
\Phi_{ji}\in \left(H^0(\mathcal{O}(-a_i+b_j+2))\otimes_\mathbb{C} Hom_\mathbb{C}(M_i, N_j)\right)^H 
\end{cases}$$

\noindent for all $1\leq i\leq m$ and $1\leq j\leq n$.
\end{theorem}
\begin{proof}
\noindent The proof is similar to the proofs of Proposition 3.5 and Theorem 3.6 in \cite{BKP}. Let us consider actions of $\beta^{-1}(H)$ on $H^0(\mathcal{O}(b-a+2))$ and on $Hom_{\mathbb{C}}(N, M)$. In both cases $-I_2$ acts via multiplication of $-1$ on $H^0(\mathcal{O}(b-a+2))$ and on $Hom_{\mathbb{C}}(N, M)$. Thus the action of $-I_2$ on $H^0(\mathcal{O}(b-a+2))\otimes_\mathbb{C} Hom_\mathbb{C}(N, M)$ is trivial and the action of $\beta^{-1}(H)$ reduces to an action of $H$.\\

\noindent The rest of the proof uses techniques similar to the ones in Theorem \ref{maintheorem1} and Remark \ref{directsum} and the fact that a $\beta^{-1}(H)$-equivariant co-Higgs field $\Phi$ on a direct sum of $\beta^{-1}(H)$-equivariant bundles $E\oplus F$ decomposes into diagonal and cross-diagonal $\beta^{-1}(H)$-equivariant components. In particular, on the spaces $$H^0(\mathcal{O}(a_i-b_j+2)), Hom_\mathbb{C}(N_j, M_i), H^0(\mathcal{O}(-a_i+b_j+2)), Hom_\mathbb{C}(M_i, N_j)$$ the element $-I_2$ acts via the multiplication by $-1$ and on respective tensor products $H^0(\mathcal{O}(a_i-b_j+2))\otimes_\mathbb{C} Hom_\mathbb{C}(N_j, M_i)$ and $H^0(\mathcal{O}(-a_i+b_j+2))\otimes_\mathbb{C} Hom_\mathbb{C}(M_i, N_j)$ it acts trivially. Hence follows the proof of the theorem.
\end{proof}
\printbibliography

\end{document}